\documentclass[11pt]{article} 
\usepackage{amsmath,amsfonts,amssymb}
\usepackage{url}
\usepackage{amsmath}

\usepackage{graphicx}
\usepackage{capt-of}
\usepackage{booktabs}
\usepackage{varwidth}
\usepackage{multirow}

\usepackage{epsfig}
\usepackage{epstopdf}
\usepackage{color}
\usepackage{ textcomp }
\usepackage{mathabx}
\usepackage{indentfirst}
\def\argmin{\mathop{\rm argmin}\nolimits}

\def\tto{\;{\lower 1pt \hbox{$\rightarrow$}}\kern -10pt
\hbox{\raise 2pt \hbox{$\rightarrow$}}\;}

\def\ra{\rangle}
\def\la{\langle}

\def\B{I\!\!B}
\def\h{\hfill\Box}
\def\R{\Bbb R}
\def\N{\Bbb N}
\def\ox{\bar{x}}
\def\ow{\bar{w}}

\def\oz{\bar{z}}

\def\ou{\bar{u}}
\def\oq{\bar{q}}

\def\dom{\mbox{\rm dom}\,}

\def\cone{\mbox{\rm cone}\,}

\def\i{\mbox{\rm int}\,}

\def\h{\hfill\square}

\def\O{\Omega}

\newcounter{lk}

\begin{document}

\begin{center}
\vspace*{0.3in} {\bf SMOOTHING ALGORITHMS FOR COMPUTING THE PROJECTION ONTO A MINKOWSKI SUM OF CONVEX SETS}\\[2ex]
Xiaolong Qin\footnote{Institute of Fundamental and Frontier Sciences, University of Electronic Science and Technology of China, Chengdu 611731, China (email: qxlxajh@163.com)}, Nguyen Thai An\footnote{Institute of Research and Development, Duy Tan University, Danang, Vietnam (email: thaian2784@gmail.com)}
\end{center}
{\bf Abstract.} In this paper, the problem of computing the projection, and therefore the minimum distance, from a point onto a Minkowski sum of general convex sets is studied. Our approach is based on the minimum norm duality theorem originally stated by Nirenberg and the Nesterov smoothing techniques. It is shown that projection points onto a Minkowski sum of sets can be represented as the sum of points on constituent sets so that, at these points, all of the sets share the same normal vector which is the negative of the dual solution. The proposed NESMINO algorithm improves the theoretical bound on number of iterations from $O(\frac{1}{\epsilon})$ by Gilbert [SIAM J. Contr., vol. 4, pp. 61--80, 1966] to $O\left(\frac{1}{\sqrt{\epsilon}}\ln(\frac{1}{\epsilon})\right)$, where $\epsilon$ is the desired accuracy for the objective function. Moreover, the algorithm also provides points on each component sets such that their sum is equal to the projection point.\\[1ex] 

\newtheorem{Theorem}{Theorem}[section]
\newtheorem{Proposition}[Theorem]{Proposition}
\newtheorem{Remark}[Theorem]{Remark}
\newtheorem{Lemma}[Theorem]{Lemma}
\newtheorem{Corollary}[Theorem]{Corollary}
\newtheorem{Definition}[Theorem]{Definition}
\newtheorem{Example}[Theorem]{Example}
\renewcommand{\theequation}{\thesection.\arabic{equation}}
\normalsize

\section{Introduction}
Let $A$ and $B$ be two subsets in $\R^n$. Recall that  the {\em Minkowski sum} of these sets is defined by
$$A+B:=\{a+b :\; a\in A, b\in B\}.$$
The case of more than two sets is defined in the same way by induction. Note that if all the sets $A_i$ for $i=1, \ldots, m$ are convex, then every  linear combination of these sets, $\sum_{i=1}^m \lambda_iA_i$ with $\lambda_i\in \R$ for $i=1, \ldots, m$, is also convex. The Euclidean distance function associated with a subset $Q$ is defined by
$$d(x; Q):=\inf\{\|q-x\|:\; q\in Q\},$$
where $\|\cdot\|$ is the Euclidean norm. The optimization problem we are concerned with in this paper is the following {\em minimum norm problem}
\begin{equation}
d\left(0; \sum_{i=1}^p T_i(\Omega_i)\right):= \min\left\{\|x\|: \; x\in \sum_{i=1}^p T_i(\Omega_i)\right\},
\label{main_prob}
\end{equation}
where $\O_i$, for $i=1,\ldots, p$, are nonempty convex compact sets in $\R^n$ and $T_i: \; \R^m \to \R^n$ are affine mappings satisfying $T_i(x)=A_ix+a_i$, where $A_i$, for $i=1, \ldots, p$, are $n\times m$ matrices and $a_i$, for $i=1, \ldots, p$, are given points in $\R^n$.
Since $\sum_{i=1}^p T_i(\Omega_i)$ is closed and convex, and the norm under consideration is Euclidean, \eqref{main_prob} has a unique solution, which is the projection from the origin onto $\sum_{i=1}^p T_i(\Omega_i)$. We denote this solution by $x^*$ throughout this paper.


We assume in problem \eqref{main_prob} that each constituent set $\O_i$ is simple enough so that the corresponding projection operator $P_{\O_i}$ is easy to compute. It is worth noting that there hasn't been an algorithm for finding the Minkowski sum of general convex sets except for the cases of balls and polytopes. Moreover, in general, the projection onto a Minkowski sum of sets cannot be represented as the sum of the projections onto constituent sets.

Minimum norm problems for the case of polytopes have been well studied in the literature from both theoretical and numerical point of view; see e.g., \cite{MDM,W,Z} and the references therein. The most suitable algorithm for solving \eqref{main_prob} is perhaps the one suggested by Gilbert \cite{G}. The original Gilbert's algorithm was devised for solving the minimum norm problem associated with just one convex compact set. The algorithm does not require the explicit projection operator of the given set. Instead, it requires in each step the computation of the support point of the set along with a certain direction. By observation that for a given direction, support point of a
Minkowski sum of sets can be represented in term of support points of constituent sets, Gilbert's algorithm thus can be applied for general case of \eqref{main_prob}. Following \cite{G}, Gilbert's algorithm is a descent method that generates a sequence $\{z_k\}$ satisfying $\|z_k\|$ converges downward to $\|x^*\|$  within $O(\frac{1}{k})$ iterations.

Another effective algorithm for distance computation between two convex objects is the GJK algorithm proposed by Gilbert, Johnson and Keerthi \cite{GJK} and its enhancing versions \cite{GF,Cam,VDB}. The original GJK algorithm was just restricted to compute the distance between objects which can be approximately represented as convex polytopes. In order to reduce the error of the polytope approximations in finding the minimum distance, Gilbert and Fo \cite{GF} modified the original GJK to handle general convex objects. The new modified algorithm is based on Gilbert's algorithm and has the same bound on number of iterations. It has been observed by many authors that, the algorithm often makes rapid movements toward the solution at its starting iterations, however on many problems, the algorithm turns to stuck when it approaches the final solution; see \cite{Keerthi,Martin}. 


When deal with problem \eqref{main_prob}, we are interested in the following questions:

Question 1: Is it possible to characterize points on each of constituent sets $\O_i$ so that the sum of their images under corresponding affine mappings is equal to $x^*$.

Question 2: Is there an alternative algorithm that improves the theoretical complexity bound of Gilbert's algorithm for solving \eqref{main_prob}?

In this work, we first use the minimum norm duality theorem originally stated by Nirenberg \cite{nb} to establish the Fenchel dual problem for \eqref{main_prob}. From this duality result, we show that each projection onto a Minkowski sum of sets can be represented as the sum of points on constituent sets so that at these points, all the sets share the same normal vector. For numerically solving the problem, we utilize the smoothing technique developed by Nesterov \cite{n,nbook}. To this end, we first approximate the dual objective function by a smooth and strongly convex function and solve this dual problem via a fast gradient scheme. After that we show how an approximate solution for the primal problem, i.e., an approximation for the projection of the origin, can be reconstructed from the dual iterative sequence. Our algorithm improves the theoretical bound on number of iterations from $O(\frac{1}{\epsilon})$ of the Gilbert algorithm to $O\left(\frac{1}{\sqrt{\epsilon}}\ln(\frac{1}{\epsilon})\right)$. Moreover, the algorithm also provides elements on each constituent sets such that the sum of their images under corresponding linear mappings is equal to the projection $x^*$.

The rest of the paper is organized as follows. In section 2, we provide tools of convex analysis that are widely used in the sequel. The Nesterov's smoothing technique and fast gradient method are recalled in section 3. In section 4, we state some duality results concerning the minimum norm problems and give the answer for Question 1. Section 5  is devoted to an overview of Gilbert's algorithm. Section 6 is the main part of the paper devoted to develop a smoothing algorithms for solving \eqref{main_prob}. Some illustrative examples are provided in section 7.   

\section{Tools of Convex Analysis}
\label{sec:pre}
In $n$-dimensional Euclidean space $\R^n$, we use $\la \cdot, \cdot \ra$ to denote the inner product, and $\|\cdot\|$ to denote the associated Euclidean norm. 
An extended real-valued function $f:\R^n \to \R\cup\{+\infty\}$ is said to be convex if  
$f((1-\lambda) x + \lambda y) \leq (1-\lambda) f(x) + \lambda f(y)$, for all $x, y \in \R^n$ and $\lambda \in (0,1)$. We say that $f$ is {\em strongly convex} with modulus $\gamma$ if $f-\frac{\gamma}{2}\|\cdot\|^2$ is a convex function. Let $Q$ be a subset of $\R^n$, the support function of $Q$ is defined by 
\begin{equation}
\sigma_{Q}(u):=\sup\{\la u, x\ra:\; x\in Q\}, \; u\in \R^n.
\label{supportfunction}
\end{equation}
It follows directly from the definition that $\sigma_{Q}$ is positive homogeneous and subadditive. We denote by $S_Q(u)$ the  solution set of \eqref{supportfunction}. The set-valued mapping $S_Q: \R^n \rightrightarrows \R^n$ is called the {\em support point mapping} of $Q$. If $Q$ is compact, then $\sigma_{Q}(u)$ is finite and $S_Q(u) \neq \emptyset$. Moreover, an element $s_Q(u)\in S_Q(u)$ is a point in $Q$ that is farthest in the direction $u$. Thus $s_{Q}(u)$ satisfies
$$s_{Q}(u)\in Q \mbox{ and } \sigma_{Q}(u)=\la u, s_{Q}(u) \ra.$$

In order to study minimum norm problem in which the Euclidean distance is replaced by distances generated by different norms, we consider a more general setting. Let $F$ be a closed, bounded and convex set of $\R^n$ that contains the origin as an interior point. The {\em minimal time function} associated with the dynamic set $F$ and the target set $Q$ is defined by 
\begin{equation}
T_F(x; Q):=\inf\{t\geq 0 :\; (x+tF)\cap Q\neq \emptyset\}.
\label{mtf}
\end{equation}
The minimal time function \eqref{mtf} can be expressed as 
\begin{equation} \label{mnt}
T_F(x; Q)=\inf\{\rho_F(\omega-x):\;\omega\in Q\},
\end{equation}
where $\rho_F(x):=\inf\{t\geq 0: \; x\in tF\}$ is the Minkowski function associated with $F$. Moreover, $T_F(\cdot, Q)$ is convex if and only if $Q$ is convex; see \cite{MN2010}. We denote by
$$\Pi_F(x; Q):=\{q\in Q: \rho_F(q-x)=T_F(x; Q)\}$$
the set of {\em generalized projection} from $x$ to $Q$.

Note that, if $F$ is the closed unit ball generated by some norm $\vvvert \cdot \vvvert$ on $\R^n$, then we have $\rho_F=\vvvert \cdot \vvvert$, $\sigma_F=\vvvert \cdot \vvvert_*$ and $T_F(\cdot,Q)$ reduces to the ordinary distance function
$$d(x; Q)=\inf\{\vvvert\omega-x\vvvert: \; \omega \in Q\}, \quad x\in \R^n.$$
The set $\Pi_F(x; Q)$ in this case is denoted by $\Pi(x; Q):=\{q \in Q :\; d(x; Q)=\vvvert q-x\vvvert\}.$ When $\vvvert \cdot \vvvert$ is Euclidean norm, we simply use the notation $P_{\O}(x)$ instead. If $\O$ is a nonempty closed convex set, then the Euclidean projection $P_{\O}(x)$ is a singleton for every $x\in \R^n$.  

The following results whose proof can be found in \cite{hul} allow us to represent support functions of general sets in term of the support functions of one or more simpler sets. 
\begin{Lemma}\label{lem:3}Consider the support function \eqref{supportfunction}. Let $\O$, $\O_1$,
	$\O_2$ be subsets of $\R^m$ and $T: \R^m \to \R^n$ satisfying $T(x)=Ax+a$ be an affine transformation, where $A$ is an $n\times m$ matrix and $a\in \R^n$. The following assertions hold: \\
	{\rm (i)} $\sigma_\O=\sigma_{\mbox{\rm cl\;}\O}=\sigma_{\mbox{\rm co\;}\O}=\sigma_{\mbox{$\overline{\rm co}\;$} \O}$.\\
	{\rm (ii)} $\sigma_{\O_1+\O_2}(u)=\sigma_{\O_1}(u)+\sigma_{\O_2}(u)$ and
	$\sigma_{\O_1 - \O_2}(u)=\sigma_{\O_1}(u) + \sigma_{\O_2}(-u)$, for all  $u\in \R^m$.\\
	{\rm (iii)} $\sigma_{T(\O)}(v)=\sigma_\O(A^\top v) +\la v, a\ra$, for all  $v\in \R^n$.
\end{Lemma}

From Lemma \ref{lem:3}, we have the following properties for support point mappings. 

\begin{Lemma}\label{lem2}
	Let $\O$, $\O_1$,
	$\O_2$ be convex compact subsets of $\R^m$ and $T: \R^m \to \R^n$ be an affine transformation satisfying $T(x)=Ax+a$ , where $A$ is an $n\times m$ matrix and $a\in \R^n$. The following assertions hold: \\
	{\rm (i)} $S_{\O_1+\O_2}(u)=S_{\O_1}(u) + S_{\O_2}(u),\; \mbox{ for all } u\in \R^m.$\\
	{\rm (ii)} $S_{T(\O)}(v)=T\left(S_{\O}(A^\top v)\right)=A\left(S_{\O}(A^\top v)\right) +a, \; \mbox{ for all } v\in \R^n.$\\
	{\rm (iii)} If suppose further that $\O$ is a strictly convex set, then $S_{\O}(u)$ is a singleton for any $u\in \R^m\setminus \{0\}$. 
\end{Lemma}
{\bf Proof. } (i) The assumption on the compactness ensures the nonemptyness of involving support points sets. Let any support point $\ow \in S_{\O_1+\O_2}(u)$. This means that $\ow \in \O_1+\O_2$ and $\la u, \ow\ra =\sigma_{\O_1+\O_2}(u)$. There exists $\ow_1 \in \O_1$ and $\ow_2\in \O_2$ such that $\ow=\ow_1+\ow_2 $. Employing Lemma \ref{lem:3}(ii), we have   
	\begin{equation}\label{sum5}
	\la u, \ow_1\ra +\la u,\ow_2\ra =\sigma_{\O_1}(u)+\sigma_{\O_2}(u).
	\end{equation}
	From the definition of support functions, $\la u, \ow_1\ra \leq \sigma_{\O_1}(u)$ and $\la u, \ow_2\ra \leq \sigma_{\O_2}(u)$. Therefore, equality \eqref{sum5} holds if and ony if $\la u, \ow_1\ra = \sigma_{\O_1}(u)$ and $\la u, \ow_2\ra = \sigma_{\O_2}(u)$. Thus, $\ow \in S_{\O_1}(u) + S_{\O_2}(u)$ and we have justified the $"\subset"$ inclusion in (i). The converse implication is straightforward. 
	
	Now let $\ow\in \O$ such that $A\ow + a \in S_{T(\O)}(u)$. By Lemma \ref{lem:3}(iii), we have $\la v, A\ow + a\ra = \sigma_{T(\O)}(v)=\sigma_\O(A^\top v) +\la v, a\ra$. This is equivalent to $\la A^\top v, \ow \ra=\sigma_\O(A^\top v)$. Therefore, $\ow\in S_{\O}(A^\top v)$ and thus $A\ow + a \in A\left(S_{\O}(A^\top v)\right) +a$. The converse implication of (ii) is proved similarly. 
	
	For (iii), suppose that there exist $\ow_1, \ow_2\in \O$ with $\ow\neq \ow_2$ such that $\la u, \ow_1\ra =\la u, \ow_2\ra =\sigma_{\O}(u)$. It then follows from the properties of support function that 
	$$\la u, \dfrac{\ow_1+\ow_2}{2} \ra =  \dfrac{1}{2}\la u, \ow_1\ra + \dfrac{1}{2}\la u, \ow_2\ra=\dfrac{1}{2}\sigma_{\O}(u) +\dfrac{1}{2}\sigma_{\O}(u)=\sigma_{\O}(u).$$
	Since $\O$ is strictly convex, $\ow:=\frac{\ow_1+\ow_2}{2} \in \i(\O)$. Take $\epsilon>0$ small enough such that $\B(\ow; \epsilon) \subset \O$. Then the element $\widehat{w}:=\ow+ \frac{\epsilon}{2}\frac{u}{\|u\|} \in \B(\ow; \epsilon) \subset \O$. Moreover, since $u\neq 0$, we have
	$$\la u, \widehat{w}\ra = \la u, \ow\ra + \frac{\epsilon}{2}\|u\| > \la u, \ow\ra =\sigma_\O(u).$$
	This is a contradiction. The proof is complete. $\h$

The \emph{Fenchel conjugate} of a convex function $f: \R^n \to \R\cup \{+\infty\}$ is defined by
\begin{equation*}
f^*(v):=\sup\{\la v,x\ra - f(x):\; x\in \R^n\}, \; v\in \R^n.
\end{equation*}
If $f$ is proper and lower semicontinuous, then $f^*: \R^n \to \R\cup \{+\infty\}$ is also a proper, lower semicontinuous convex function. From the definition, support function $\sigma_{Q}$ is the Fenchel conjugate of the {\em indicator function} $\delta_{Q}$ of $Q$ which is defined by $\delta_{Q}(x)=0$ if $x\in Q$ and $\delta_Q(x)=+\infty$ otherwise. 

The polar of a subset $E\subset \R^n$ is the set $E^\circ=\{u\in \R^n:\; \sigma_E(u) \leq 1\}$. When $E$ is the closed unit ball of a norm $\vvvert \cdot\vvvert$, then $E^\circ$ is the closed unit ball of the corresponding dual norm $\vvvert\cdot\vvvert_{*}$.  Some basis properties of the polar set are collected in the following result whose proof can be found in \cite[Proposition 1.23]{Tuy}.
\begin{Proposition}The following assertions hold:\\\label{prop_polar}
	{\rm (i)} For any subset $E$, the polar $E^\circ$ is a closed convex set containing the origin and $E\subset E^{\circ\circ}$;\\
	{\rm (ii)} $0 \in \i(E)$ if and only if $E^\circ$ is bounded;\\
	{\rm (iii)} $E = E^{\circ\circ}$ if $E$ is closed convex and contains the origin;\\
	{\rm (iv)} If $E$ is closed convex and contains the origin, then $\rho_E=\sigma_{E^{\circ}}$ and $\left(\rho_E\right)^*=\delta_{E^{\circ}}$.
\end{Proposition}
Thus, if $F$ is a closed convex and bounded set with $0\in \i(F)$ then $F^\circ$ is also a closed convex and bounded set with $0\in \i(F^\circ)$. Moreover, from the subadditive property, $\rho_F=\sigma_{F^\circ}$ is a Lipschitz function with modulus $\|F^\circ\|:=\sup\{\|x\|: \; x\in F^\circ\}$.

Let us recall below the Fenchel duality theorem which plays an important role in the sequel. We denote the set of points where a function $g: \R^n\to \R \cup\{+\infty\}$ is finite and continuous by $\mbox{\rm cont} g$. 
\begin{Theorem}{\rm (See \cite[Theorem 3.3.5]{BL})} Given functions $f: \R^m \to \R \cup\{+\infty\}$ and $g: \R^n\to \R \cup\{+\infty\}$, and a linear mapping $A: \R^m \to \R^n$, the weak duality inequality
	\begin{align*}
	\inf \limits_{x\in \R^m}\{f(x) + g(Ax)\} \geq \sup \limits_{u\in \R^n} \{ -f^*(A^*u)-g^*(-u)\}
	\end{align*}
	holds. If furthermore $f$ and $g$ are convex and satisfy the following condition 
	\begin{equation}
	A(\dom f) \cap \mbox{\rm cont} g \neq \emptyset,
	\label{QC}
	\end{equation}
	then the equality holds and the supremum is attained if it is finite. \label{Fthm} 
\end{Theorem}

\section{Nesterov's Smoothing Technique and Fast Gradient Method}\label{sec:Nes}
In a celebrated work \cite{n83}, Nesterov introduced a fast firtst-order method for solving convex smooth problems in which the objective functions have Lipschitz continuous gradients. In contrast to the complexity bound of $O(1/\epsilon)$ possessed by the classical  gradient descent method, Nesterov's method gives a complexity bound of $O(1/\sqrt{\epsilon})$, where $\epsilon$ is the desired accuracy for the objective function.

When the problem under consideration is nonsmooth in which the objective function has an explicit max-structure as follows
\begin{equation}\label{orig_func}
f(u):=\max\{\la Au, x\ra -\phi(x) :\; x\in Q\},\; u\in \R^n,
\end{equation}
where $A$ is an $m\times n$ matrix and $\phi$ is a continuous convex function on a compact set $Q$ of $\R^m$, in order to overcome the complexity bound $O(\frac{1}{\epsilon^2})$ of the subgradient method, Nesterov \cite{n} made use of the special structure of $f$ to approximate it by a function with Lipschitz continuous gradient and then applied a fast gradient method to minimize the smooth approximation. With this combination, we can solve the original non-smooth problem up to accuracy $\epsilon$ in $O(\frac{1}{\epsilon})$ iterations. To this end, let $d$ be a continuous strongly convex function on $Q$. Let $\mu$ be a positive number called a \emph{smooth parameter}. Define
\begin{equation}\label{smoothing_approx}
f_\mu(u):=\max\{\la Au, x\ra-\phi(x)-\mu d(x) :\; x\in Q\}.
\end{equation}
Since $d(x)$ is strongly convex, problem \eqref{smoothing_approx} has a unique solution. The following statement is a simplified version of \cite[Theorem~1]{n}.
\begin{Theorem}{\rm (See \cite[Theorem 1]{n})}\label{s1_I} The function $f_\mu$ in {\rm (\ref{smoothing_approx})} is well defined and continuously differentiable on $\R^n$. The gradient of the function is
	$$\nabla f_\mu(u)=A^\top x_\mu(u),$$
	where $x_\mu(u)$ is the unique element of $Q$ such that the maximum in {\rm (\ref{smoothing_approx})} is attained. Moreover, $\nabla f_\mu$ is a Lipschitz function with the Lipschitz constant $\ell_\mu=\dfrac{1}{\mu\sigma_1}\|A\|^2,$
	and
	$$f_\mu(u)\leq f(u)\leq f_\mu(u)+\mu D\quad \forall u\in \R^n,$$
	where $D:=\max\{d(x) :\; x\in Q\}$.
\end{Theorem}

For the reader's convenience, we conclude this section with a presentation of the simplest optimal method for minimizing smooth strongly convex functions; see \cite{nbook} and the references therein. Let $g: \R^n \to \R$ be strongly convex with parameter $\gamma>0$ and its gradient be Lipschitz continuous with constant $L>\gamma$. Consider   problem $g^*=\inf\left\{ g(u): u\in \R^n\right\}
$ and denote by $u^*$ its unique optimal solution.

\medskip
{\begin{center}
		\begin{tabular}{| l |}
			\hline
			\qquad \qquad {\bf Fast Gradient Method}\\
			\hline
			{\small INITIALIZE}: $\gamma$, $v_0=u_0 \in \R^n$.\\
			Set $k=0$.\\
			{\bf Repeat the following}\\
			\qquad Set $u_{k+1}:=v_k - \frac{1}{L}\nabla g(v_k)$\\
			\qquad Set $v_{k+1}:=u_{k+1} + \frac{\sqrt{L} -\sqrt{\gamma}}{\sqrt{L} + \sqrt{\gamma}}\left(u_{k+1}-u_k\right)$\\
			\qquad Set $k:=k+1$\\
			{\bf Until a stopping criterion is satisfied.}\\
			\hline
		\end{tabular}
\end{center}}

\medskip
By taking into account \cite[Theorem 2.2.3]{nbook},   $\{u_k\}_{k=0}^\infty$ satisfies
\begin{align}
\label{est2g}
g(u_k)-g^*& \leq \left(g(u_0) -g^* + \frac{\gamma}{2}\|u_0-u^*\|^2 \right)\left(1-\sqrt{\frac{2}{L}}\right)^k \notag\\
& \leq \left(g(u_0) -g^* + \frac{\gamma}{2}\|u_0-u^*\|^2 \right)e^{-k\sqrt{\frac{\gamma}{L_\mu}}} \notag \\
& \leq 2\left(g(u_0) -g^* \right)e^{-k\sqrt{\frac{\gamma}{L}}},
\end{align}
while the last inequality is a consequence of  \cite[Theorem 2.2.3]{nbook}.

Since $g$ is a differentiable strongly convex function and $u^*$ is its unique minimizer on $\R^n$, we have $\nabla g\left(u^*\right)=0$. Using \cite[Theorem 2.1.5]{nbook}, we find
\begin{equation}\label{est3g}
\dfrac{1}{2L}\|\nabla g(u_k)\|^2 \leq g(u_k) - g^* \stackrel{\eqref{est2g}}{\leq} 2\left(g(u_0) -g^* \right)e^{-k\sqrt{\frac{\gamma}{L}}}.
\end{equation}


\section{Duality for Minimum Norm Problems}
\label{sec:dual}
In this section, we are in a position to give some duality results concerning minimum norm problem \eqref{main_prob}. Let us first recall the duality theorem originally stated by Nirenberg \cite{nb}.
\begin{Theorem}{\rm (Minimum norm duality theorem)} \label{thmdual} Given $\ox\in \R^n$ and let $d(\cdot; \O)$ be the distance  function to a nonempty closed convex set $\O$ associated with some norm $\vvvert \cdot \vvvert$ on $\R^n$. Then
	$$d(\ox; \O)=\max\{\la u, \ox\ra -\sigma_\O(u): \; \vvvert u \vvvert_{*} \leq 1\},$$
	where the maximum on the right is achieved at some $\ou$. Moreover, if $\ow\in\Pi(\ox; \O)$, then $-\ou$ is aligned with $\ow -\ox$, i.e., $\la -\ou, \ow -\ox \ra =\vvvert \ou \vvvert_* . \vvvert \ow-\ox \vvvert$.
\end{Theorem}

\begin{figure}[!ht]
	\centering
	\vspace{-6cm}
	\includegraphics[width=15cm]{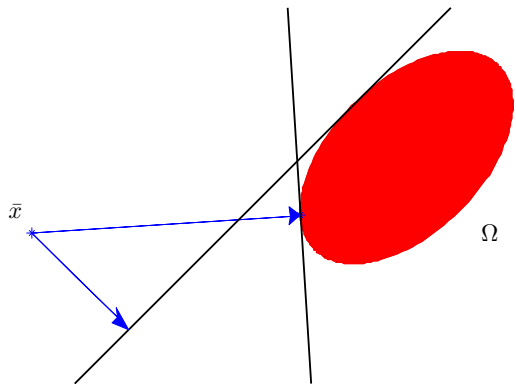}
	\vspace{-7.5cm}
	\caption{\small An illustration of the minimum norm duality theorem}\label{MNT}
\end{figure}

According to this theorem, the minimum distance from a point to a convex set is equal to the maximum of the distance from the point to hyperplanes separating the point and the set; see Figure \ref{MNT}. A standard proof of this theorem can be found in \cite[p. 136]{Lube}. We also refer the readers to the recent paper \cite{Dax} for more types of minimum norm duality theorems concerning the width and the length of symmetrical convex bodies.

\begin{Lemma}\label{nonempty}
	Let $Q$ be a nonempty closed subset of $\R^n$. Then the generalized projection $\Pi_F(x; Q)$ is nonempty for any $x\in \R^n$.
\end{Lemma}

{\bf Proof.} From the assumption that $F$ is a closed bounded and convex set that contains the origin as an interior point, $0\leq T_F(x; Q) <+\infty$ for all $x\in \R^n$ and the following number exists 
	$$R=\sup\{r:\; \B(0; r)\subset F^\circ\}<+\infty.$$
	Then we have $\rho_F(x)=\sigma_{F^\circ}(x) \geq R\|x\|$ for all $x\in \R^n$.
	Fix $x\in \R^n$. For each $n\in \N$, from \eqref{mnt} there exists $w_n\in Q$, such that 
	
	\begin{equation} \label{est10}
	T_F(x; Q) \leq \rho_F(w_n-x) < T_F(x; Q) +\frac{1}{n}.
	\end{equation}
	It follows from \eqref{est10} and triangle inequality that 
	that 
	$$R\|w_n\| \leq  R\left(\|w_n-x\| +\|x\| \right)\leq \rho_F(w_n-x)+R\|x\|\leq  T_F(x;Q) +1 +R\|x\|$$
	for all $n$. Thus the sequence $\{w_n\}$ is bounded. We can take a subsequence $\{w_{k_n}\}$ that converges to a point $\ow\in Q$ due to the closedness of $Q$. By taking the limit both sides of \eqref{est10} and using the continuity of $T_F(\cdot; Q)$ and the Minkowski function, we can conclude that $\ow \in \Pi_F(x; Q)$. $\h$ 

\medskip
Theorem \ref{thmdual} is in fact a direct consequence of the Fenchel duality theorem which is used to prove the following extension for minimal time functions.

\begin{Theorem}\label{thm_dual}The generalized distance $T_F(0; A(\O))$ from the origin $0_{\R^n}$ to the image $A(\O)$ of a nonempty closed convex set $\O \subset \R^m$ under a linear mapping $A: \R^m \to \R^n$ can be computed by
	$$T_F(0; A(\O)):=\inf\{\rho_F(Aw):\; w\in \O\} =\max\{-\sigma_\O(-A^\top u): \;u\in F^\circ\},$$
	where the maximum on the right is achieved at some $\ou\in F^\circ$.  If $A\ow \in \Pi_F\left(0; A(\O)\right)$ is a projection from the origin to $A(\O)$, then
	\begin{equation*}
	\la A\ow, \ou\ra = \sigma_{F^\circ}(A\ow) =-\sigma_\O(-A^\top \ou).
	\end{equation*}
\end{Theorem}
{\bf Proof.} Applying Theorem \ref{Fthm} for $g=\rho_F$ and $f=\delta_\O$, the following qualification condition holds
	$$A\mbox{dom} f \cap \mbox{cont}(g) = A(\O)\cap \R^n =A(\O)\neq \emptyset,$$
	where $\mbox{cont}(g)=\R^n$ is due to the fact that $\rho_F$ is a continuous function on $\R^n$. It follows that
	\begin{align*}
	T_F(0; A(\O))&=\min\{\delta_{\O}(x) + \rho_F(Ax): \;x\in \R^n\}\\
	&=\sup\{-\left(\delta_\O\right)^*\left(A^\top u\right) - \left(\rho_F\right)^*(-u):\; u\in \R^n\}\\
	&=\sup\{-\sigma_{\O}\left(A^\top u\right) - \delta_{F^\circ}(-u):\; u\in \R^n\}\\
	&=\sup\{-\sigma_{\O}\left(-A^\top u\right) - \delta_{F^\circ}(u):\; u\in \R^n\}\\
	&=\sup\{-\sigma_{\O}\left(-A^\top u\right):\; u\in F^\circ\},
	\end{align*}
	and the supremum is attained because $T_F(0; A(\O))$ is finite. If the supremum on the right is achieved at some $\ou\in F^\circ$ and the  infimum on the left is achieved at some $\ow\in \O$, then
	$$\sigma_{F^\circ}\left(A\ow\right)=\rho_F(A\ow) = T_F(0, A(\O))=\max\{-\sigma_\O\left(-A^\top u\right):\; u\in F^\circ\}=-\sigma_\O(-A^\top \ou).$$
	Since $\ow\in \O$, we also have
	$$\la -A^\top \ou, \ow\ra \leq \sigma_\O\left(-A^\top \ou\right)=-\sigma_{F^\circ}(A\ow).$$
	This implies that  $\la A^\top \ou, \ow\ra \geq \sigma_{F^\circ}(A\ow).$ On the other hand,
	$\sigma_{F^\circ}(A\ow) \geq \la A\ow, \ou \ra = \la A^\top \ou, \ow \ra$, because $\ou\in F^\circ$. Thus,
	$\la A\ow, \ou\ra = \sigma_{F^\circ}(A\ow)$. This completes the proof. $\h$

Note that, given a closed set $\O$, the set $A(\O)$ need not to be closed and therefore, we can not use the $\min$ to replace the $\inf$ in the primal problem in Theorem \ref{thm_dual}.
\begin{Proposition}\label{dual_prop} Let $Q$ be a nonempty,  closed convex subset of $\R^n$. The following holds
	\begin{equation} \label{primdual}
	T_F(0; Q):=\min\{\rho_F(q):\; q\in Q\} =\max\{-\sigma_Q(-u): \;u\in F^\circ\}.
	\end{equation}
	If the maximum on the right is achieved at $\ou\in F^\circ$ and the infimum on the left is attained at $\oq\in Q$, then
	\begin{equation}
	\la \oq, \ou\ra = \sigma_{F^\circ}(\oq) =-\sigma_Q(-\ou).
	\label{parallel}
	\end{equation}
	If  $F=\B$ is the Euclidean closed unit ball , then the projection $\oq$ exists uniquely and 
	$$d(0; Q):=\min\{\|q\|:\; q\in Q\} =\max\{-\sigma_Q(-u): \;u\in \B\}.$$
	If suppose further that $0\notin Q$, then $\frac{\oq}{\|\oq\|}$ is the unique solution of the dual problem.
\end{Proposition}
{\bf Proof.} The first assertion is a direct consequence of Theorem \ref{thm_dual} with $\O=Q$ and $A$ is the identity mapping of $\R^n$. Note that, by Lemma \ref{nonempty}, the infimum is also attained here. When $F$ is the Euclidean ball, the minimal time function reduces to the Euclidean distance function and therefore the projection $\oq=P_Q(0)$ exists uniquely. If $0\notin Q$, then $\oq\neq0$. Moreover, we have $\la -\oq, x-\oq\ra \leq 0$ for all $x\in Q$. This implies,
	$$\left \la -\frac{\oq}{\|\oq\|}, x \right \ra \leq -\|\oq\|, \mbox{ for all } x\in Q.$$
	Hence $\sigma_Q\left(-\frac{\oq}{\|\oq\|}\right)\leq -\|\oq\|=-d(0; Q)$. This means that $\frac{\oq}{\|\oq\|}$ is a solution of the following dual problem
	$$d(0; Q)=\max\{-\sigma_Q(-u): \;u\in \B\}.$$
	From \eqref{parallel}, any dual solution $\ou$ must satisfy $\ou\in S_{F^\circ}(\oq)$. Since $F=\B$, we have $F^\circ=\B$ is a strictly convex set. Thus, by Lemma \ref{lem2}(iii),  $\ou=\frac{\oq}{\|\oq\|}$ is the unique solution of dual problem. The proof is now complete. $\h$

\begin{figure}[!ht]
	\centering
	\vspace{-5.5cm}
	\includegraphics[width=13cm]{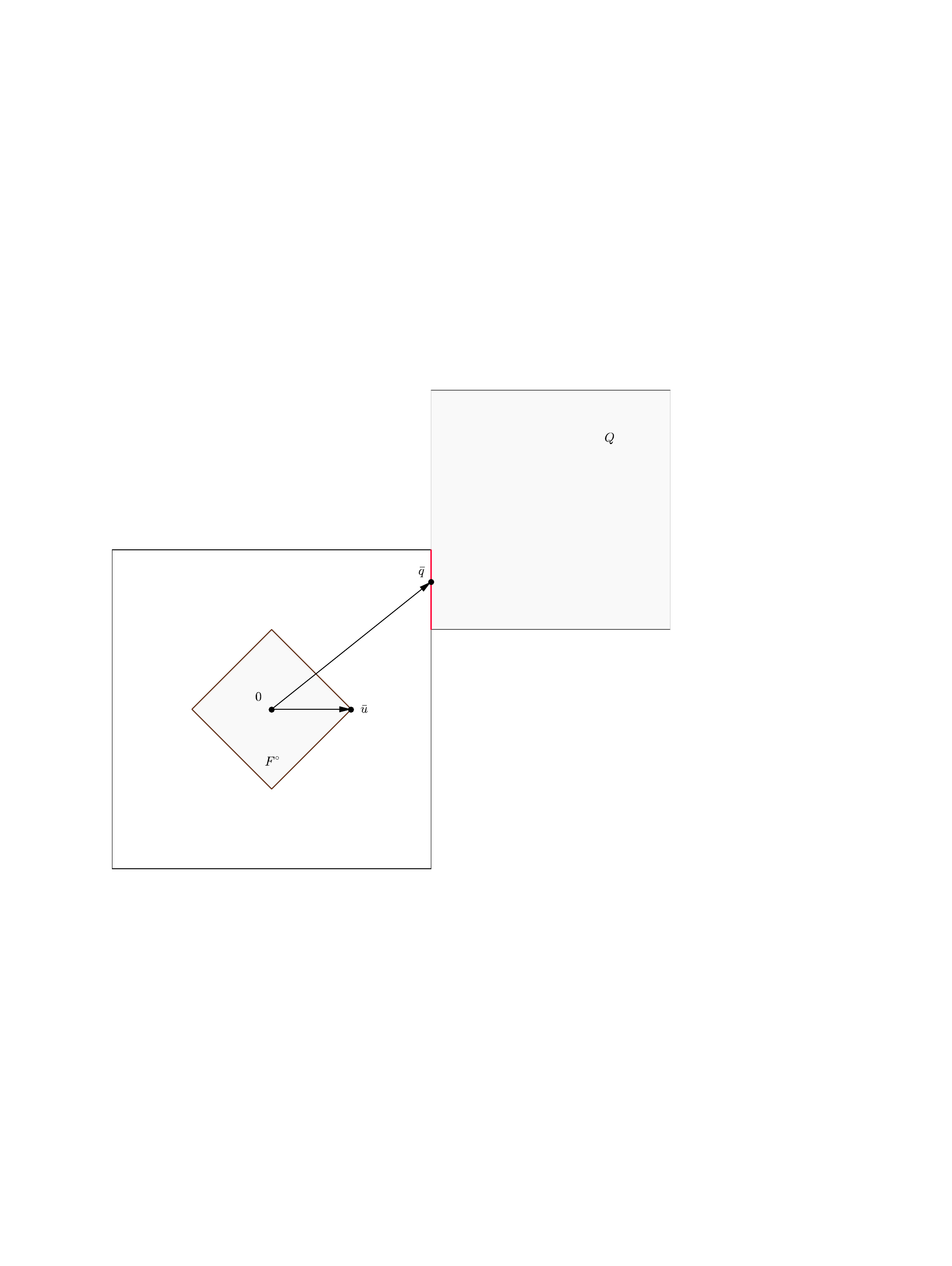}
	\vspace{-5.5cm}
	\caption{\small A minimum norm problem with non-Euclidean distance.}
	\label{FnonEu}
\end{figure}

\medskip
From \eqref{parallel}, for any primal-dual pair $(\oq, \ou)$, we have the following relationship
\begin{equation}\label{par2}
\ou \in S_{F^\circ}(\oq) \;\;\; \mbox{   and   }\;\;\; \oq \in S_{Q}(-\ou).
\end{equation}
This observation seems to be useful from numerical point of view in the sense that if a dual solution $\ou$ is found exactly, then a primal solution $\oq$ can be obtained by taking a support point in $S_Q(-\ou)$.  However, for a general convex set $Q$, the set $S_Q(-\ou) $ might contain more than one point and there might be some points in this set which is not a desired primal solution. Thus, the above task is possible when $S_Q(-\ou)$ is a singleton.    

When the distance function under consideration is non-Euclidean, the primal problem may have infinitely many solutions and we may not recover a dual solution from a primal one $\oq$ by setting $\frac{\oq}{\|\oq\|}$ as in the Euclidean case.

\begin{Example}{\rm 
		In $\R^2$, consider the problem of finding the projection onto the set $Q=\{x\in \R^2:\; 2\leq x_1\leq 5 \mbox{ and } 1\leq x_2\leq 4\}$ in which the   
		distance function generated by the $\ell_\infty$-norm. In this case, $F=\{x\in \R^2:\; \max\{|x_1|,|x_2|\}\leq 1\}$ and 
		$F^\circ=\{x\in \R^2:\; |x_1|+|x_2|\leq 1\}$ and we have $T_F(0; Q)=2$. The primal problem in \eqref{primdual} has the solution set 
		$\Pi_F(0; Q)=\{x\in \R^2:\; x_1=2 \mbox{  and } 1\leq x_2\leq 2\}$ and the corresponding dual problem has a unique solution $\ou=(1,0)$. We can see that, for any primal solution $\oq$, the element $\frac{\oq}{\|\oq\|} \neq \ou $. Thus, $\frac{\oq}{\|\oq\|}$ is not a dual solution; see Figure \ref{FnonEu}.     
	}
\end{Example}

We now give a sufficient condition for the uniqueness of solution of primal and dual problems in \eqref{primdual}. We recall the following definition from \cite{nars}. The set $F$ is said to be {\em normally smooth} if and only if for every boundary point $\ox$ of $F$, the normal cone of $F$ at $\ox$ defined by $N(\ox; F):=\{u\in \R^n:\; \la u, x-\ox\ra \leq 0, \;\forall x\in F \}$ is generated exactly by one vector. That means, there exists $a_{\ox}\in \R^n$ such that $N(\ox; F)=\cone\{a_{\ox}\}$. From \cite[Proposition 3.3]{nars}, we have that $F$ is normally smooth if and only if its polar $F^\circ$ is
strictly convex.

\begin{Proposition} We have the following: \\
	{\rm (i)} If $Q$ is a nonempty closed and strictly convex set of $\R^n$, then the generalized projection set $\Pi_F(x; Q)$ is a singleton for all $x\in \R^n$.\\
	{\rm (ii)} If $F$ is normally smooth, then the dual problem in \eqref{primdual} has a unique solution.
\end{Proposition}
{\bf Proof.} (i) It follows from the definitions of minimal time function and generalized projection that $T_F(x; Q)=T_F(0; Q-\{x\})$ and $\Pi_F(x; Q)=x+\Pi_F(0; Q-\{x\})$. It suffices to prove that $\Pi_F(0; Q-\{x\})\neq \emptyset$.  Let $\ou$ is a dual solution in \eqref{primdual}. Since $Q$ is nonempty and closed, the set $\Pi_F(0; Q-\{x\})$ is nonempty by Lemma \ref{nonempty}. Suppose that $\Pi_F(0; Q-\{x\})$ contains two distinct elements $q_1-x\neq q_2-x$. Then, by relation \eqref{par2}, both $q_1$ and $q_2$ belong to the set $S_Q(-\ou)$. This is a contradiction to Lemma \ref{lem2} by the strictly convexity of $Q$ and justifies (i). The proof of (ii) is similar by using the strictly convexity of $F^\circ$. $\h$  

\medskip
Minkowski sum of two closed sets is not necessarily closed. For example, for $Q_1=\{x\in \R^2:\; x_2\geq e^{x_1}\}$ and $Q_2=\{x\in \R^2:\; x_2=0\}$, the sum  $$Q_1+Q_2=\{x\in \R^2:\; x_2>0 \}$$ is an open set. In what follows, in order to ensure the existence of support point for the Minkowski sum, we assume that all component sets are compact.   

\begin{figure}[!ht]
	\centering
	\vspace{-6.4cm}
	\includegraphics[width=14cm]{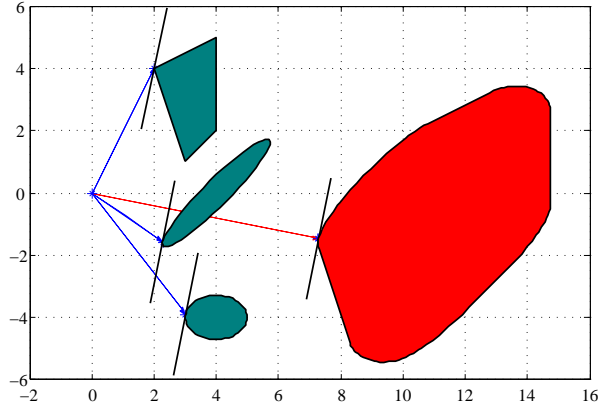}
	\vspace{-6cm}
	\caption{\small The Minkowski sum of  a polytope and two ellipses is approximately plotted by the red set. The projection of the origin onto the red is the sum of points on three constituent sets such that, at these points, three sets have a normal vector in common.}
	\label{F1}
\end{figure}

We now show that \eqref{par2} can allow us to characterize points on each constituent sets in the Minkowski sum so that their sum is equal to the projection point. The answer for Question 1 in the introduction is stated in the following results; see Figure \ref{F1} for an illustration.

\begin{Corollary}\label{Corf}
	Let $\{Q_i\}_{i=1}^p$ be a finite collection of nonempty convex compact sets in $\R^n$. It holds that
	$$T_F\left(0; \sum_{i=1}^p Q_i\right) = -\min\left\{\sum_{i=1}^p\sigma_{Q_i}(-u): \; u\in F^\circ\right\}.$$
	Moreover, if the minimum on the right hand side is attained at $\ou\in F^\circ$, then any generalized projection $\oq$ of the origin onto the set $\sum_{i=1}^p Q_i$ satisfies 
	$$\oq\in S_{Q_1}(-\ou)+\ldots+S_{Q_p}(-\ou).$$
	Thus, the projection $\oq$ is the sum of points on component sets such that at these points all the sets have the same normal vector $-\ou$.
	If  $F=\B$ is the Euclidean closed unit ball , then the projection $\oq$ exists uniquely and 
	$$d\left(0; \sum_{i=1}^p Q_i\right) = -\min\left\{\sum_{i=1}^p\sigma_{Q_i}(-u): \; u\in \B\right\}.$$
	If in addition, $0\notin \sum_{i=1}^p Q_i$ then $\frac{\oq}{\|\oq\|}$ is the unique solution of the dual problem and  we have $$\oq \in S_{Q_1}(-\oq)+\ldots+S_{Q_p}(-\oq).$$
\end{Corollary}
{\bf Proof.} Let $Q:=\sum_{i=1}^p Q_i$. Using Lemma \ref{lem:3} and Lemma \ref{lem2}, we have
	$$\sigma_Q(-u)=\sigma_{Q_1}(-u)+\ldots+\sigma_{Q_p}(-u)  \;  \mbox{   and   }\;   S_Q(-u)=S_{Q_1}(-u)+\ldots+ S_{Q_p}(-u).$$
	Note that, the support point mapping $S_{Q}(u)$ does not depend on the magnitude of $u$, using Proposition \ref{dual_prop} and relation \eqref{par2}, we clarify the desired conclusion easily. $\h$

\medskip
The problem of finding a pair of closest points, and therefore the Euclidean distance, between two given convex compact sets $\mathcal{P}$ and $\mathcal{Q}$ can be reduced to the minimum problem associated with the Minkowski sum $\mathcal{Q}-\mathcal{P}$ by observing that $d(\mathcal{P}, \mathcal{Q})=d(0, \mathcal{Q}-\mathcal{P})$.  A note  here is that although there may be several pairs of closest points, the latter problem always has a unique solution which is the projection from $0$ onto $\mathcal{Q}-\mathcal{P}$.   By noting that $\sigma_{-\mathcal{P}}(-u)=\sigma_{\mathcal{P}}(u)$ and $S_{-\mathcal{P}}(-u)=-S_{\mathcal{P}}(u)$, we have the following result; see Figure \ref{F4}.

\begin{Corollary}
	Let $\{Q_i\}_{i=1}^p$ and $\{P_j\}_{j=1}^l$ be two finite collection of nonempty convex compact sets in $\R^n$ and let $\mathcal{P}=\sum_{j=1}^l P_j, \mathcal{Q}= \sum_{i=1}^p Q_i$.  It holds that
	$$d\left(\mathcal{P}, \mathcal{Q}\right) = -\min\left\{\sum_{i=1}^p\sigma_{Q_i}(-u) + \sum_{j=1}^l\sigma_{P_j}(u): \; u\in \B\right\}.$$
	Moreover, if $\oq$  is the projection of the origin onto $Q:=\mathcal{Q}-\mathcal{P}$, if $(\bar{a}, \bar{b})$ is a pair of closest points of  $\mathcal{Q}$ and $\mathcal{P}$, then $\oq=\bar{a} - \bar{b}$ and 
	$$\bar{a}\in S_{Q_1}(-\oq)+\ldots+S_{Q_p}(-\oq) \;\; \mbox{    and    }\;\; \bar{b}\in S_{P_1}(\oq)+\ldots+S_{P_\ell}(\oq).$$
	Thus, $\bar{a}$ is the sum of points in $Q_i$ for $i=1, \ldots, p$ such that at these points all $Q_i$ have the same normal vector $-\oq$ and $\bar{b}$ is the sum of points in $P_j$ for $j=1, \ldots, \ell$ such that at these points all $P_j$ have the same normal vector $\oq$.
\end{Corollary}

\section{The Gilbert Algorithm}
\label{sec:G}
We now give an overview and clarify how the Gilbert algorithm can be applied for solving \eqref{main_prob}. Let us define the function $g: \R^n \times Q \to \R$ by
\begin{equation}
g_Q(z, x):=\sigma_Q(z) - \la z, x\ra,
\label{function_g}
\end{equation}
where $Q=\sum_{i=1}^p T_i(\O_i)$. From the definition, $g_Q(-z, z)\geq 0$ for all $z\in Q$. A point  $z\in Q$ is the solution of \eqref{main_prob} if and only if  $\la -z, x-z\ra \leq 0$ for all $x\in Q$. This amounts to saying that $g_Q(-z, z)=0$.

\begin{Lemma}\label{segment}
	If two points $z$ and $\oz$ satisfy $\|z\|^2 -\la z, \oz \ra >0$, then there is a point $\tilde{z}$ in the line segment $\mbox{co}\{z, \oz\}$ such that $\|\tilde z\|< \|z\|$.
\end{Lemma}
{\bf Proof.} If $\|\oz\|^2\leq\la z, \oz\ra$, then we can choose $\tilde{z}=\oz$. Consider the case $\|\oz\|^2>\la z, \oz\ra$. By combining with the assumption  $\|z\|^2 -\la z, \oz \ra >0$, we have
	$$0<\lambda^*:=\dfrac{\|z\|^2 -\la z, \oz\ra}{\|z-\oz\|^2} <1.$$
	This implies the quadratic function
	$$f(\lambda)=\|\oz -z\|^2\lambda^2 + 2\la z, \oz-z\ra \lambda +\|z\|^2$$
	attains its minimum on $[0, 1]$ at $\lambda^*$ and therefore $f(\lambda^*)=\|z+\lambda^*(\oz-z)\|^2 < f(0)=\|z\|^2.$ Thus $\tilde{z}:=z+\lambda^*(\oz-z)$ is the desired point. $\h$

\begin{figure}[!ht]
	\centering
	\vspace{-0.4cm}
	\includegraphics[width=10cm]{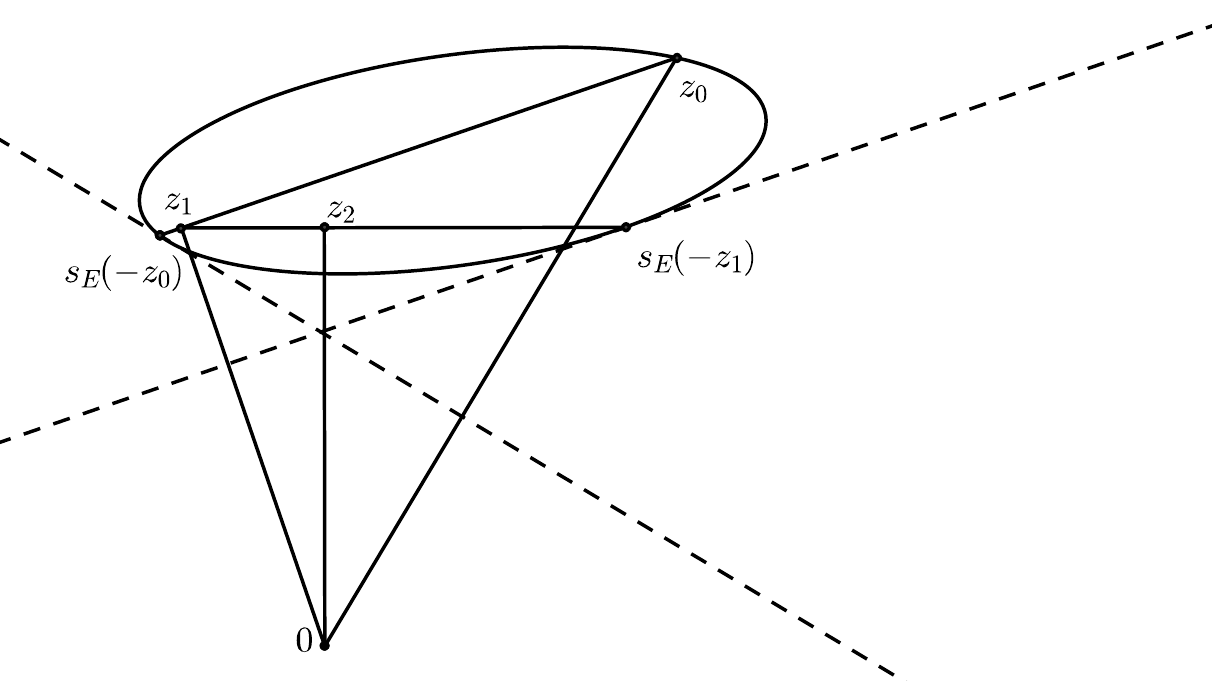}
	\caption{\small An illustration of Gilbert's algorithm.}\label{FGA}
\end{figure}

The Gilbert algorithm can be interpreted as follows. Starting from some $z\in Q$, if $g_Q(-z, z)=0$ then $z$ is the solution. If $g_Q(-z, z)>0$, then $\oz\in S_Q(-z)$ satisfies $\|z\|^2 -\la z, \oz \ra >0$. Using Lemma \ref{segment}, we   find a point $\tilde{z}$ on the line segment connecting  $z$ and $\oz$ such that  $\|\tilde{z}\|<\|z\|$.  The algorithm is outlined  as follows.

\medskip
{\begin{center}
		\begin{tabular}{| l |}
			\hline
			\qquad \qquad \qquad \qquad  \quad {\bf Gilbert's Algorithm}\\
			\hline
			0. Initialization step: Take arbitrary point $z_0\in Q$.\\
			1. If $g_Q(-z,z)=0$, then return $z=x^*$ is the solution\\
			\qquad \qquad \qquad \qquad   \quad else, set $\oz\in S_Q(-z)$.\\
			2. Compute $\tilde{z}\in \mbox{co}\{z, \oz\}$ which has minimum norm, set $z=\tilde{z}$ \\
			and go back to step 1.\\
			\hline
		\end{tabular}
\end{center}}

\medskip
Figure \ref{FGA} illustrates some iterations of Gilbert's algorithm for finding closest point to an ellipse in two dimension. Lemma \ref{segment} also suggests an effective way to find   $\tilde{z}$ in step 3. We have $\tilde{z}:=z+\lambda^*(\oz-z)$, where
$$\lambda^*=\begin{cases}
1,&  \mbox{if } \|\oz\|^2\leq\la z, \oz\ra, \\
\dfrac{\|z\|^2 -\la z, \oz\ra}{\|z-\oz\|^2},& \mbox{otherwise}.
\end{cases}$$
To implement the algorithm, it remains to show how to compute a supporting point for $Q=\sum_{i=1}^p T_i(\O_i)$. Fortunately, this can be done by using Lemma \ref{lem2}.

Gilbert showed that, if   $\{z_k\}_{k=1}^\infty$ generated by the algorithm does not stop with $z=x^*$ at step 1 within a finite number of iterations, then $z_k\to x^*$ asymptotically.  According to \cite[Theorem 3]{G}, we have
\begin{equation}\label{boundG}
\|z_k\| - \|x^*\| \leq \frac{C_1}{k}  \quad \mbox{and } \quad \|z_k - x^*\| \leq \frac{C_2}{\sqrt{k}},
\end{equation}
where $C_1$ and $C_2$ are some positive constants. From the above estimates, in order to find an $\epsilon$ - approximate solution, i.e., a point $z$ such that $\|z\| - \|x^*\|\leq \epsilon$, we need to perform the algorithm in $O(\frac{1}{\epsilon})$ iterations. Gilbert also showed the bounds $\eqref{boundG}$ are sharp in the sense that within a constant multiplicative factor it is imposible to obtain  bounds on $\|z_k\| - \|x^*\|$ and $\|z_k - x^*\| $ which approach zero more rapidly than those given \eqref{boundG}; see \cite[Example 1]{G}.

\section{Smoothing Algorithm for Minimum Norm Problems}
\label{sec:main}
Our approach for numerically solving \eqref{main_prob} is based on the minimum norm duality Theorem \ref{thm_dual} and the Nesterov smoothing technique \cite{n}. Let us first consider the function of the following type
\begin{equation*}
\sigma_{A, Q}(u)=\sup\{\la Au, x\ra: \; x\in Q\},\; u\in \R^n,
\end{equation*}
where $A$ is an $m\times n$ matrix and $Q$ is a closed bounded subset of $\R^m$. Observe that $\sigma_{A, Q}(u)$ is the composition of  a linear mapping and the support function of $Q$. As we will see, this function can be approximated by the following function
\begin{equation*}
\sigma^\mu_{A, Q}(u)=\sup \left\{\la Au, x\ra - \frac{\mu}{2}\|x\|^2: \; x\in Q \right\},\; u\in \R^n.
\end{equation*}
The following statement is a directly consequence of Theorem \ref{s1_I}. However, the approximate function as well as its gradient, in this case, has closed form that is expressed in term of the Euclidean projection. This feature makes it reliable from numerical point of view.
\begin{Proposition}\label{mainp}
	The function $\sigma^\mu_{A, Q}$ has the following explicit representation
	$$\sigma^\mu_{A, Q}(u)=\dfrac{\|Au\|^2}{2\mu} -\dfrac{\mu}{2}\big[d(\dfrac{Au}{\mu}; Q)\big]^2$$
	and is continuous differentiable on $\R^n$ with its gradient given by
	$$\nabla \sigma^\mu_{A, Q}(u)=A^\top P_Q\left(\dfrac{Au}{\mu}\right).$$
	The gradient $\nabla \sigma^\mu_{A, Q}$ is a Lipschitz function with constant $\ell_\mu=\dfrac{1}{\mu}\|A\|^2 $. Moreover,
	\begin{equation}
	\sigma^\mu_{A, Q}(u)\leq \sigma_{A, Q}(u)\leq \sigma^\mu_{A, Q}(u)+\dfrac{\mu}{2} \|Q\|^2\; \mbox{\rm for all }u\in \R^n,
	\label{est0}
	\end{equation}
	where $\|Q\|:=\sup\{\|q\| :\; q\in Q\}.$
\end{Proposition}
{\bf Proof.} We have
	\begin{align*}
	\sigma^\mu_{A, Q}(u)&=\sup\left \{ \la Au, x\ra -\dfrac{\mu}{2}\|x\|^2 :\; x\in Q \right\}\\
	&=\sup\left \{-\dfrac{\mu}{2}\big(\|x\|^2-\dfrac{2}{\mu}\la Au, x\ra\big) :\; x\in Q \right\}\\
	&=-\dfrac{\mu}{2}\inf \left\{\|x-\dfrac{Au}{\mu}\|^2-\dfrac{\|Au\|^2}{\mu^2} :\; x\in Q \right\}\\
	&=\dfrac{\|Au\|^2}{2\mu} - \dfrac{\mu}{2}\inf \left\{\|x-\dfrac{Au}{\mu}\|^2 :\; x\in Q\right\}\\
	&=\dfrac{\|Au\|^2}{2\mu} -\dfrac{\mu}{2}\left[d\left(\dfrac{Au}{\mu}; Q\right)\right]^2.
	\end{align*}
	Since $\psi(x):=[d(x; Q)]^2$ is a differentiable function satisfying $\nabla \psi(x)=2[x-P_Q(x)]$ for all $x\in \R^m$, we find from the chain rule that
	\begin{align*}
	\nabla \sigma^\mu_{A, Q}(u)&=\dfrac{1}{\mu}A^\top(Au)  -\dfrac{\mu}{2}\left[\dfrac{2}{\mu}A^\top\left(\dfrac{Au}{\mu} - P_Q(\dfrac{Au}{\mu})\right)\right]\\
	&=A^\top P_Q(\dfrac{Ax}{\mu}).
	\end{align*}
	From the property of the projection mapping onto convex sets and Cauchy-Schwarz inequality, we find,  for any $u, v\in \R^n$, that 
	\begin{align*}
	\|\nabla \sigma^\mu_{A, Q}(u)  - \nabla \sigma^\mu_{A, Q}(v)\|^2
	&=\|A^\top P_Q(\dfrac{Au}{\mu}) - A^\top P_Q(\dfrac{Av}{\mu})\|^2\\
	&\leq \|A\|^2 \|P_Q(\dfrac{Au}{\mu}) - P_Q(\dfrac{Av}{\mu}) \|^2\\
	&\leq \|A\|^2 \left\la\dfrac{Au-Av}{\mu}, P_Q(\dfrac{Au}{\mu}) - P_Q(\dfrac{Av}{\mu})\right\ra\\
	&=\dfrac{\|A\|^2}{\mu}\left \la u-v, A^\top P_Q(\dfrac{Au}{\mu}) - A^\top P_Q(\dfrac{Av}{\mu})\right \ra\\
	&=\dfrac{\|A\|^2}{\mu}\la u-v, \nabla \sigma^\mu_{A, Q}(u) - \nabla \sigma^\mu_{A, Q}(v)\ra\\
	&\leq\dfrac{\|A\|^2}{\mu}\|u-v\| \|\nabla \sigma^\mu_{A, Q}(u) - \nabla \sigma^\mu_{A, Q}(v)\|.
	\end{align*}
	This implies that
	$$\|\nabla \sigma^\mu_{A, Q}(u) - \nabla \sigma^\mu_{A, Q}(v)\| \leq  \dfrac{\|A\|^2}{\mu}\|u-v\|.$$
	The lower and upper bounds in~\eqref{est0} follow from
	\begin{align*}
	\la Au, x\ra -\dfrac{\mu}{2}\|x\|^2 &\leq \la Au, x\ra \leq \la Au, x\ra -\dfrac{\mu}{2}\|x\|^2+ \dfrac{\mu}{2}\sup\left\{\|q\|^2 :\; q\in Q\right\},
	\end{align*}
	for all $x\in Q$. The proof is now complete. $\h$

\medskip
From Proposition \ref{dual_prop}, we have the duality result below
$$d\left(0; \sum_{i=1}^p T_i(\O_i)\right) = -\min\left\{\sum_{i=1}^p\sigma_{\O_i}(-A_i^\top u) - \la u, \sum_{i=1}^p a_i\ra:\; u\in \B\right\}.$$

We now make use of the strong convexity of the squared Euclidean norm to state another dual problem for \eqref{main_prob} in which the dual objective function is strongly convex. 
\begin{Proposition} The following duality result holds
	\begin{equation}
	d^2\left(0; \sum_{i=1}^p T_i(\O_i)\right)= - \min\{\sum_{i=1}^p\sigma_{\O_i}\left({A_i}^\top u\right) + \la u, \sum_{i=1}^p a_i\ra+ \frac{1}{4}\|u\|^2:\; u\in \R^n\}.
	\label{dual_Eu}
	\end{equation}
\end{Proposition}
{\bf Proof.} Observe that the Fenchel dual of the function $\|\cdot\|^2$ is   $\frac{1}{4}\|\cdot\|^2$. Applying Theorem \ref{Fthm} again, for $Q:=\sum_{i=1}^p T_i(\O_i)$, we have
	\begin{align*}
	\left[d\left(0; Q\right)\right]^2 &= \inf\left\{\|x\|^2:\; x\in Q          \right\}\\
	&=\max\{-\left(\delta_Q\right)^*\left( u\right) - \left(\|\cdot\|^2 \right)^*(-u): \; u\in \R^n\}\\
	&=\max\{-\sigma_{Q}\left( u\right) - \frac{1}{4}\|-u\|^2: \; u\in \R^n\}\\
	&=-\min\{\sigma_{Q}\left(u\right) + \frac{1}{4}\|u\|^2: \; u\in \R^n \}.
	\end{align*}
	The result now follows directly from Lemma \ref{lem:3}. $\h$

\medskip
In order to solve   minimum norm problem \eqref{main_prob}, we solve   dual problem \eqref{dual_Eu} by approximating the dual objective function by a smooth and strongly convex function with Lipschitz continuous gradient and then apply  a fast gradient scheme to this smooth one.

Let us define the dual objective function by
$$f(u):=\sum_{i=1}^p \sigma_{\O_i}\left({A_i}^\top u\right) + \la u,\sum_{i=1}^p a_i \ra+ \frac{1}{4}\|u\|^2, \quad u\in \R^n.$$
The following result is a direct consequence of Proposition \ref{mainp}.
\begin{Proposition}
	The function $f(u)$ has the following smooth approximation
	$$f_\mu(u):=\sum_{i=1}^p \left(\dfrac{\|A_i^\top u\|^2}{2\mu} - \dfrac{\mu}{2}[d(\dfrac{A_i^\top u}{\mu}; \O_i)]^2\right) + \la u,\sum_{i=1}^p a_i \ra+ \frac{1}{4}\|u\|^2, \quad u\in \R^n.$$
	Moreover, $f_\mu$ is a strongly convex function with modulus $\gamma=2$ and its gradient is given by
	\begin{equation}
	\nabla f_\mu(u)=\sum_{i=1}^p A_i P_{\O_i}\left(\dfrac{A_i^\top u}{\mu}\right) +\sum_{i=1}^p a_i+ \frac{1}{2}u.
	\label{nabla}
	\end{equation}
	The Lipschitz constant of $\nabla f_\mu$ is
	\begin{equation}
	L_\mu := \frac{\sum_{i=1}^p \|A_i\|^2}{\mu} +\frac{1}{2}.
	\label{Lmu}
	\end{equation}
	Moreover, we have the following estimate
	\begin{equation}
	f_\mu(u) \leq f(u) \leq f_\mu(u) +\mu D_f,
	\label{est3}
	\end{equation}
	where $D_f:=\dfrac{1}{2}\sum_{i=1}^p  \|\O_i\|^2<\infty$.
\end{Proposition}

We now apply the Nesterov fast gradient method introduced in Section 3 for minimizing the smooth and strongly convex function $f_\mu$. We will show how to recover an approximately optimal solution for primal problem   \eqref{main_prob} from the dual iterative sequence. The {\bf NE}sterov {\bf S}moothing algorithm for {\bf MI}nimum {\bf NO}rm Problem (NESMINO) is outlined as follows:

\smallskip
{\begin{center}
		\begin{tabular}{| l |}
			\hline
			\qquad \qquad \qquad \quad \qquad {\bf NESMINO}\\
			\hline
			{\small INITIALIZE}:   $\O_i, A_i, a_i$ for $i=1, \ldots, p$ and  $v_0$, $u_0$, $\mu$. \\
			Set $k=0$.\\
			{\bf Repeat the following}\\
			\qquad Compute $\nabla f_\mu(v_k)$ using \eqref{nabla}\\
			\qquad Compute $L_{\mu}$ using \eqref{Lmu}\\
			\qquad Set $u_{k+1}:=v_k - \frac{1}{L_\mu}\nabla f_\mu(v_k)$\\
			\qquad Set $v_{k+1}:=u_{k+1} + \frac{\sqrt{L_\mu} -\sqrt{2}}{\sqrt{L_\mu} + \sqrt{2}}\left(u_{k+1}-u_k\right)$\\
			\qquad Set $k:=k+1$\\
			{\bf Until a stopping criterion is satisfied.}\\
			\hline
		\end{tabular}
\end{center}}

\smallskip
We denote by $u^*_\mu$ the unique minimizer of $f_\mu$ on $\R^n$. We also denote by $u^*$ a minimizer of $f$ and by $f^*:=f(u^*)=\inf_{x\in \R^n} f(x)$ its optimal value on $\R^n$. From the duality result \eqref{dual_Eu}, we have
$$f^*=-\left[d\left(0; \sum_{i=1}^p T_i(\O_i)\right)\right]^2.$$
We say that $x\in \sum_{i=1}^pT_i(\O_i)$ is an  $\epsilon$-approximate solution of problem \eqref{main_prob} if it satisfies
$$\|x \| - \|x^*\| \leq \epsilon.$$

\begin{Theorem}
	Let $\{u_k\}_{k=1}^\infty$ be the sequence generated by {\rm NESMINO} algorithm. Then  the sequence $\{y_k\}_{k=1}^\infty$ defined by $$y_k:=\sum_{i=1}^p \left[A_iP_{\O_i}\left(\frac{A_i^\top u_k}{\mu}\right)+a_i\right]$$ converges to an  $\epsilon$-approximate solution of  minimum norm problem \eqref{main_prob} within  $k=O\left(\frac{1}{\sqrt{\epsilon}}\ln\left(\frac{1}{\epsilon}\right)\right)$ iterations.
\end{Theorem}
{\bf Proof.} Using \eqref{est2g}, we find that  $\{u_k\}_{k=0}^\infty$ satisfies
	\begin{align}
	\label{est2}
	f_\mu(u_k)-f_\mu^*& \leq 2\left(f_\mu(u_0) -f_\mu^* \right)e^{-k\sqrt{\frac{\gamma}{L_\mu}}}.
	\end{align}
	From  $f_\mu(u_0) \leq f(u_0)$ and the following estimate
	$$f_\mu^* =f_\mu(u^*_\mu) \geq f(u^*_\mu) - \mu D_f  \geq f(u^*) - \mu D_f=f^* - \mu D_f,$$
	we have
	\begin{equation}
	f_\mu(u_0) - f^*_\mu \leq  f(u_0) - f^* +\mu D_f.
	\label{est1}
	\end{equation}
	Moreover, since $f_\mu(u_k) - f^*_\mu \geq f(u_k) -\mu D_f - f^*$, we find from \eqref{est2} and \eqref{est1} that
	\begin{align}
	f(u_k) - f^*& \leq \mu D_f + f_\mu(u_k) - f^*_\mu \notag \\
	&\leq \mu D_f + 2\left(f(u_0) - f^* + \mu D_f\right)e^{-k\sqrt{\frac{\gamma}{L_\mu}}}, \mbox{ for all } k\geq 0. \label{est4}
	\end{align}
	Since $f_\mu$ is a differentiable strongly convex function and $u^*_\mu$ is its unique minimizer on $\R^n$, we have $\nabla f_\mu\left(u^*_\mu\right)=0$. It follows from \eqref{est3g} that 
	$$\dfrac{1}{2L_\mu}\|\nabla f_\mu(u_k)\|^2 \leq f_\mu(u_k) - f^*_\mu \stackrel{\eqref{est2}}{\leq} 2\left(f_\mu(u_0) -f_\mu^* \right)e^{-k\sqrt{\frac{\gamma}{L_\mu}}}.$$
	This implies
	\begin{equation}
	\|\nabla f_\mu(u_k)\|^2 \leq 4L_\mu(f_\mu(u_0) -f_\mu^*)e^{-k\sqrt{\frac{\gamma}{L_\mu}}} \stackrel{\eqref{est1}}{\leq} 4L_\mu(f(u_0) - f^* +\mu D_f)e^{-k\sqrt{\frac{\gamma}{L_\mu}}} .
	\label{est5}
	\end{equation}
	For each $k$ and for each $i\in \{1, \ldots, m\}$, let $x_k^i$ be the unique solution to the problem
	$$\sigma_{\mu, \O_i}\left(A_i^\top u_k\right):= \sup\left\{\la A_i^\top u_k, x\ra -\frac{\mu}{2}\|x\|^2: x\in \O_i \right\}.$$
	We have
	\begin{align*}
	\sup\left\{\la A_i^\top u_k, x\ra -\frac{\mu}{2}\|x\|^2: x\in \O_i \right\}&=\sup\left\{\frac{\|A_i^\top u_k\|^2}{2\mu} -\frac{\mu}{2}\left \|x- \frac{A_i^\top u_k}{\mu} \right\|^2  : x\in \O_i \right\}\\
	&=\frac{\|A_i^\top u_k\|^2}{2\mu} - \left[d\left(\frac{A_i^\top u_k}{\mu}, \O_i\right) \right]^2.
	\end{align*}
	Hence $x_k^i=P_{\O_i}\left(\frac{A_i^\top u_k}{\mu} \right)$.
	For each $k$, set
	$$d_k:=\left\|  \sum \limits_{i=1}^p \left(A_ix_k^i+a_i\right) \right\|^2 - \left[ d\left(0, \sum \limits_{i=1}^p T_i(\O_i) \right)\right]^2.$$
	Observe  $y_k:= \sum \limits_{i=1}^p \left(A_ix_k^i+a_i\right) \in \sum \limits_{i=1}^p T_i(\O_i)$.  From the property of the projection onto convex sets, we have
	\begin{align*}
	\|y_k -x^*\|^2 &= \|y_k\|^2 - \|x^*\|^2 +2\left\la -x^*, y_k-x^*\right\ra \leq \|y_k\|^2 - \|x^*\|^2 =d_k.
	\end{align*}
	This implies that   $\{y_k\}$ converges to $x^*$ whenever $d_k \to 0$ as $k\to \infty$. Moreover, we have
	\begin{equation}
	2\|x^*\|\left(\|y_k\|-\|x^*\|\right) \leq \left(\|y_k\| + \|x^*\|\right)\left(\|y_k\|-\|x^*\|\right) =\|y_k\|^2 -\|x^*\|^2=d_k.
	\label{estf}
	\end{equation}
	We have the following
	\begin{align*}
	d_k&=\|  \sum \limits_{i=1}^p (A_ix_k^i+a_i) \|^2 +f^*\\
	&=\left\|  \sum \limits_{i=1}^p \left(A_ix_k^i+a_i\right) \right\|^2+f_\mu(u_k) +f^* - f_\mu(u_k)\\
	&=\| \sum \limits_{i=1}^p (A_ix_k^i+a_i) \|^2+\sum \limits_{i=1}^p  [\la A_i^\top u_k, x_k^i \ra -\frac{\mu}{2}\|x_k^i\|^2]+ \la u_k, \sum_{i=1}^p a_i\ra +\frac{1}{4}\|u_k\|^2 \\
	&\qquad \qquad \qquad \qquad \qquad \qquad \qquad \qquad \qquad \qquad + f^* - f_\mu(u_k)\\
	&=\|  \sum \limits_{i=1}^p (A_ix_k^i+a_i) \|^2 + \big \la u_k, \sum \limits_{i=1}^p (A_ix_k^i +a_i) \big \ra +\frac{1}{4}\|u_k\|^2  - \frac{\mu}{2}\sum \limits_{i=1}^p \|x_k^i\|^2\\
	&\qquad \qquad \qquad \qquad \qquad \qquad \qquad \qquad \qquad \qquad +f^* - f_\mu(u_k)\\
	&=\|  \sum \limits_{i=1}^p \left(A_ix_k^i+a_i\right)  +\frac{1}{2}u_k \|^2 - \frac{\mu}{2}\sum \limits_{i=1}^p \|x_k^i\|^2 +f^* - f_\mu(u_k)\\
	&=\|  \sum \limits_{i=1}^p A_iP_{\O_i}(\frac{A_i^\top u_k}{\mu} ) + \sum_{i=1}^p a_i +\frac{1}{2}u_k \|^2 - \frac{\mu}{2}\sum \limits_{i=1}^m \|x_k^i\|^2 +f^* - f_\mu(u_k)\\
	&=\| \nabla f_\mu(u_k) \|^2 - \frac{\mu}{2}\sum \limits_{i=1}^p \|x_k^i\|^2 +f^* - f_\mu(u_k).
	\end{align*}
	Observe  $|f_\mu(u_k)-f^*| \stackrel{\eqref{est3}}{\leq}  |f(u_k)-f^*|+\mu D_f$  and $\sum \limits_{i=1}^p \|x_k^i\|^2\leq 2D_f$. Taking into account \eqref{est4} and \eqref{est5}, we have
	\begin{align*}
	d_k &\leq \left\| \nabla f_\mu(u_k) \right\|^2 + |f(u_k)-f^*|+2\mu D_f\\
	&\leq  4L_\mu(f(u_0) - f^* +\mu D_f)e^{-k\sqrt{\frac{\gamma}{L_\mu}}} + \mu D_f \\
	&\qquad \qquad \qquad \qquad \qquad \qquad + 2(f(u_0) - f^* + \mu D_f)e^{-k\sqrt{\frac{\gamma}{L_\mu}}} + 2\mu D_f.\\
	&\leq  2(2L_\mu+1)\left(f(u_0) - f^* +\mu D_f\right)e^{-k\sqrt{\frac{\gamma}{L_\mu}}} + 3\mu D_f.
	\end{align*}
	Now, for a fix $\epsilon>0$, in order to achieve an $\epsilon$ - approximate solution for the primal problem, we should force each of the two terms in the above estimate less than or equal to $\frac{\epsilon}{2}$. If we choose the value of smooth parameter $\mu$ to be $\frac{\epsilon}{6D_f}$, we have $d_k\leq \epsilon$ when
	\begin{equation}
	k\geq \sqrt{\frac{L_\mu}{\gamma}}\ln\left(\dfrac{4(2L_\mu+1)\left(f(u_0)-f^*+\frac{\epsilon}{6}\right)}{\epsilon}\right),
	\label{rate}
	\end{equation}
	where $L_\mu=\frac{\sum_{i=1}^p 6\|A_i\|^2D_f}{\epsilon}+\frac{1}{2}$.
	
	Thus, from \eqref{estf}, we can find an $\epsilon$ - approximate solution for primal problem within $k=O\left(\frac{1}{\sqrt{\epsilon}}\ln\left(\frac{1}{\epsilon}\right)\right)$ iterations. The proof is complete. $\h$

\bigskip
We highlight the fact that the algorithm does not require computation of the Minkowski sum but rather only the projection onto each of the constituent sets $\O_i$. Fortunately, many useful projection operators are easy to compute. Explicit formula
for   projection operator $P_{\O}$ exists when $\O$ is a closed Euclidean ball, a closed rectangle, a hyperplane, or a half-space. Although there are no analytic solutions, fast algorithms for computing the projetion operators exist for the cases of unit simplex, the closed $\ell_1$ ball (see \cite{Duchi,Codat}), or the ellipsoids (see \cite{Dai}).

In some cases, by making use of the special structure of the support function of $\O$, we can have a suitable smoothing technique in order to avoid working with   implicit projection operator $P_{\O}$ or to employ some fast projection algorithm. We consider two important cases as follows:

{\bf The case of ellipsoids.} Consider the case of ellipsoids associated with Euclidean norm
$$E(A, c):=\left\{x\in \R^n: \; (x-c)^\top A^{-1}(x-c) \leq 1\right\},$$
where the shape matrix $A$ is positive definite and the center $c$ is some given point in $\R^n$. It is well known that the support function of this Ellipsoid is $\sigma_E(u)= \sqrt{u^\top Au} + u^\top c$  and the support point in direction $u$ is $s_E(u)=\frac{Au}{\sqrt{u^\top Au}}+c$. We can rewrite the support function as follows
$$\sigma_E(u)=\sigma_{\B}\left(A^{1/2}u\right)+u^\top c,$$
where $\B$ stands for the closed unit Euclidean ball and $A^{1/2}$ is the square root of $A$. The smooth approximation $g_{\mu}$ of   function $g = \sigma_E$ has the following explicit representation
$$g_\mu(u) = \dfrac{\|A^{1/2}u\|^2}{2\mu} -\dfrac{\mu}{2}\big[d(\dfrac{A^{1/2}u}{\mu}; \B)\big]^2+u^\top c.$$
and is differentiable on $\R^n$ with its gradient given by $\nabla g(u)=A^{1/2} P_{\B}\left(\dfrac{A^{1/2}u}{\mu}\right)+c.$
Thus, instead of projecting onto the Ellipsoid, we just need to project onto the closed unit ball.

\begin{figure}
	\hspace{-1cm}
	\begin{minipage}[b]{0.50\linewidth}
		\includegraphics[width=0.9\textwidth]{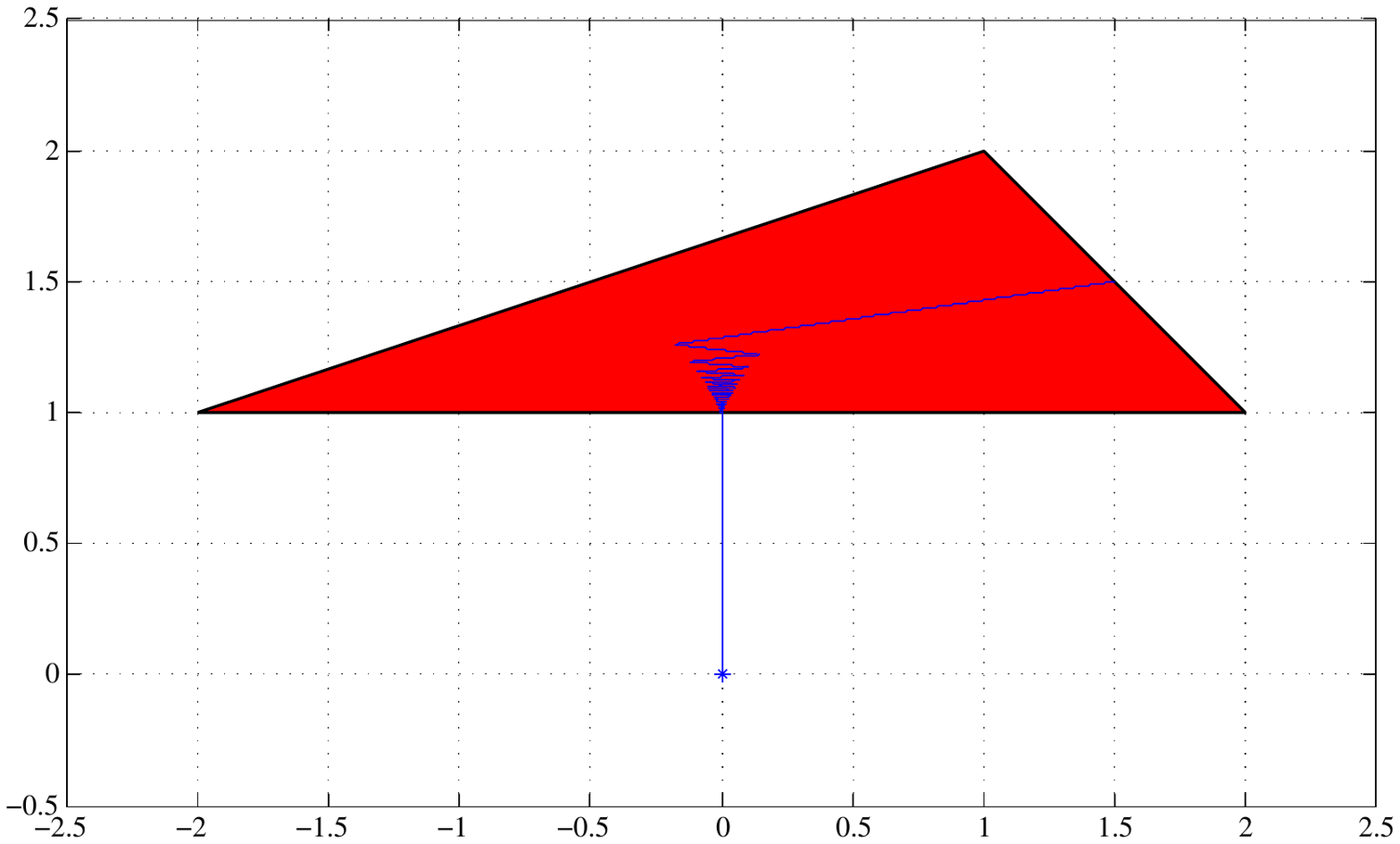}
		\vspace{-4.5cm}
	\end{minipage}%
	\begin{minipage}[b]{0.50\linewidth}
		\begin{tabular}{| c | c | c |}
			\hline
			{\rm Iteration}&{\rm NESMINO}&{\rm Gilbert}\\
			\hline
			$1$&$(0,0)$&$(1.5, 1.5)$\\
			$100$&$(0, 1)$&$(0.0090, 1.0177)$\\
			$300$&$(0, 1)$&$(0.0032, 1.0063)$\\
			$1000$&$(0, 1)$&$( 0.0010, 1.0020)$\\
			$10000$&$(0, 1)$&$( 0.0001, 1.0002)$\\
			\hline
		\end{tabular}
	\end{minipage}
	\caption{Comparison between NESMINO  algorithm and Gilbert's algorithm  for solving a minimum norm problem involving polytopes} 
	\label{F2}
\end{figure}

\medskip
{\bf The case of polytopes.} Consider the polytope $S=\mbox{conv}\{a_1, \ldots, a_m\}$ generated by $m$ point in $\R^n$. By
\cite[Theorem 32.2]{r}, we have 
$$\sigma_{S}(u)=\sup\{\la u, x\ra: \; x \in S\}=\max_{1\leq i\leq m}\la u, a_i\ra,$$ 
and the support point $s_S(u)$ of $S$ is some point $a_i$ such that $\la u, a_i\ra = \sigma_S(u)$. Observe that, for $\alpha=(\alpha_1, \ldots,\alpha_m)^\top \in \R^m$, we have
$$\max_{1\leq i\leq m}\alpha_i=\sup\{x_1\alpha_1+\ldots+x_m\alpha_m:\, x_i\geq 0, \sum_{i=1}^m x_i=1\}=\sup\left\{\la \alpha,x\ra:\, x\in \Delta_m\right\}.$$
Threfore, $\sigma_{S}(u)=\max_{1\leq i\leq m}\la u, a_i\ra=\sup\{\la Au, x\ra:\; x\in \Delta_m\}=\sigma_{\Delta_m}(Au)$, where $A=\begin{bmatrix}
a_1^\top\\
\vdots\\
a_m^\top\\
\end{bmatrix}$ and $\Delta_m$ is the unit simplex in $\R^m$. The smooth approximate function of $g=\sigma_S$ is
$g_{\mu}(u)=\dfrac{\|Au\|^2}{2\mu} -\dfrac{\mu}{2}\big[d(\dfrac{Au}{\mu}; \Delta_m)\big]^2$, with $\nabla g_\mu(u)=A^\top P_{\Delta_m}\left(\dfrac{Au}{\mu}\right).$ We thus can employ the fast and simple algorithms for computing the projection onto a unit simplex, for example in \cite{ChenYe}, instead of projection onto a polytope.

\begin{Remark}\label{Rm}
	{\rm In NESMINO algorithm, a smaller smooth parameter $\mu$ is often better because it reduces the error when approximate $f$ by $f_\mu$.
		However, a small $\mu$ implies a large value of the Lipschitz constant $L_\mu$ which in turn reduces the convergence rate by \eqref{rate}. Thus the time cost of the algorithm is expensive if we fix a value for $\mu$ ahead of time.  In practice, a sequence of smooth problems with decreasing smooth parameter $\mu$ is solved and the solution of the previous problem is used as the initial point for the next one. The algorithm stops when a preferred $\mu_*$ is attained. The optimization scheme is outlined as follows.
		
		\medskip
		{\begin{center}
				\begin{tabular}{| l |}
					\hline
					{\small INITIALIZE}: $\O_i, A_i, a_i$ for $i=1, \ldots, p$ and $w_0$, $\sigma \in (0, 1)$, $\mu_0>0$ and $\mu_*>0$.\\
					Set $k=0$.\\
					{\bf Repeat the following}\\
					\quad 1. Apply NESMINO algorithm with $\mu=\mu_k$, $u_0=v_0=w_k$ to find \\
					\quad \quad \quad \quad  \quad  \quad $w_{k+1}=\argmin_{w \in \R^n}{f_\mu(w)}$.\\
					\quad 2. Update $\mu_{k+1}:=\sigma\mu_k$ and set $k:=k+1$.\\
					{\bf Until $\mu\leq \mu_*$.}\\
					\hline
				\end{tabular}
		\end{center}}
	}
\end{Remark}

\section{Illustrative Examples}\label{sec:exp}
We now implement NESMINO and Gilbert's algorithm to solve minimum norm problem \eqref{main_prob} in a number of examples by MATLAB. We terminate the NESMINO when $\|\nabla f_\mu(u_k)\| \leq \epsilon$, for some tolerance $\epsilon>0$. In Gilbert's algorithm, we relax the stopping criterion $g_Q(-z, z)=0$ to $g_Q(-z, z)\leq\delta$, for some $\delta > 0$. The parameters described in Remark \ref{Rm} are chosen as follows:
$$\mu_0=100, \sigma=0.1, \mu_*=10^{-3}, \epsilon=10^{-3}, w_0=0,$$
and we use $\delta=10^{-4}$. All the test are implemented on a personal computer with an Intel Core i5 CPU 1.6 GHz and 4G of RAM.  Figures in this section are plotted via the Multi-Parametric Toolbox \cite{MPT} and Ellipsoidal Toolbox \cite{ET}.

Let us first give a simple example showing that when the sets involved are polytopes, the Gilbert's algorithm may have zigzag phenomenon and may become very slow as it approaches the final solution. 
\begin{Example}
	{\rm
		Consider the minimum norm problem associated with a polytope $P$ in $\R^2$ whose vertices given by the columns of the following matrix
		$$\begin{pmatrix}
		-2&2&1\\
		1&1&2\\
		\end{pmatrix}.
		$$
		The NESMINO algorithm with a fixed value $\mu=0.1$ converges to the optimal solution $x^*=(0;1)$ within nearly 100 steps. In contrast, if starting from $z_0=(\frac{3}{2}, \frac{3}{2})$, the approximate values $(z_1; z_2)$ in Gilbert's algorithm are still changing after $10^4$ iterations. In this case, as the number of iterations is increasing, the Gilbert algorithm alternately chooses the two vertices $(-2,1)$ and $(2,1)$ as support points of $P$ and converges slowly to $x^*=(0;1)$; see Figure \ref{F2}.
	}
\end{Example}
\begin{figure}[!ht]
	\centering
	\vspace{-5cm}
	\includegraphics[width=12cm]{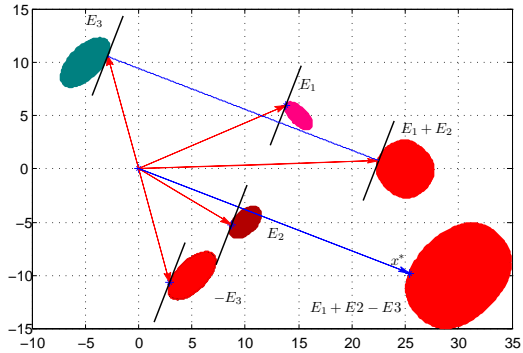}
	\vspace{-5cm}
	\caption{Minimum distance between ellipse $E_3$ and the sum $E_1+E_2$ of two other ellipses.}\label{F4}
\end{figure}

\begin{Example}
	{\rm In $\R^2$, consider the problem of computing the projection onto the sum of a polytope $P$ whose vertices given by the columns of the following matrix
		$$
		\begin{pmatrix}
		4&4&2&3\\
		2&5&4&1\\
		\end{pmatrix}$$
		and two ellipses $E_1(A_1, c_1)$, $E_2(A_2,c_2)$ with shape matrices and centers respectively given by
		$$A_1=\begin{pmatrix}
		1&0\\
		0&0.5\\
		\end{pmatrix}, c_1= \begin{pmatrix}
		4\\
		-4\\
		\end{pmatrix}
		\mbox{  and  }
		A_2=\begin{pmatrix}
		2&1\\
		1&2\\
		\end{pmatrix},
		c_2= \begin{pmatrix}
		4\\
		0\\
		\end{pmatrix}.$$
		The NESMINO algorithm yields an approximate solution $x^*=( 7.2841, -1.4787)$. The algorithm also gives 
		$x_1= (2, 4)$, $ x_2=(3.0101, -3.8995)$, $x_3 =( 2.2740, -1.5792)$ which respectively belongs to $P$, $E_1$, $E_2$ such that $x^*=x_1+x_2+x_3$. This result is depicted in Figure \ref{F1}.
	}
\end{Example}

\begin{Example}
	{\rm In this example, we apply NESMINO to find the minimum distance between a Minkowski of two ellipses $E_1(A_1, c_1), E_2(A_2, c_2)$ with shape matrices and centers respectively given by
		$$A_1=\begin{pmatrix}
		1.5&-1\\
		-1&1.5\\
		\end{pmatrix}, c_1= \begin{pmatrix}
		15\\
		5\\
		\end{pmatrix}
		\mbox{  and  }
		A_2=\begin{pmatrix}
		2&1\\
		1&2\\
		\end{pmatrix},
		c_2= \begin{pmatrix}
		10\\
		-5\\
		\end{pmatrix}$$
		and another ellipse $E_3(A_3, c_3)$ with $c_3= \begin{pmatrix}
		-5\\
		10\\
		\end{pmatrix}, A_3=\begin{pmatrix}
		5&3\\
		3&5\\
		\end{pmatrix}$. The NESMINO yields the distance  $d=  27.2347$ that is the norm of the projection $x^*= ( 25.4219, -9.7703)$ of the origin onto $E_1+E_2-E_3$. Moreover, $x^*=\bar{a} -\bar{b}$, where $\bar{a}=(22.4983, 0.8118)\in E_1+E_2$ and
		$\bar{b}=(  -2.9236, 10.5820)\in E_3$ is the pair of closest points; see Figure \ref{F4}.
	}
\end{Example}

\begin{Example} {\rm We now consider the problem of computing the projection of the origin onto a Minkowski sum of two ellipsoids $E_1(A_1, c_1)$ and $E_2(A_2, c_2)$ in high dimensions. Let $M=\left\{10^{\frac{(i-1)cond}{d-1}}:\; i=1, \ldots d\right\}$ where $d$ is the space dimension and the number $cond$ allows us to adjust the shapes (thin or fat) of the ellipsoids. For each pair $(d, cond)$, we generate 1000 problems and implement both NESMINO and Gilbert algorithm, and compute the average CPU time in seconds. In each problem, let $A_1$ and $A_2$ are $d\times d$ diagonal matrices such that the main diagonal entries of each of them is some permutation of $M$. Since we want to guarantee that $0\notin E_1+E_2$, we choose the two corresponding centers $c_1, c_2\in \R^{d\times 1}$ such that each of their entries is chosen randomly between $m$ and $11m$, where $m=\sqrt{\frac{10^{cond}}{d}}$.  The result is reported in Table \ref{Tab1}.
		
		\begin{table}[!ht]
			\caption{Performance of NESMINO algorithm and Gilbert's algorithm in solving minimum norm problems associated with ellipsoids.}\label{Tab1}
			\centering
			\renewcommand{\arraystretch}{1.25}
			\rotatebox{90}{
				\begin{tabular}{| c | c | l | l |}
					\hline
					$d$&$cond$&NESMINO&Gilbert\\
					\hline
					\multirow{4}{*}{10}&2&0.0021&0.0007\\
					&3&0.0061& 0.0006\\
					&4& 0.0192&0.0006\\
					&5& 0.0634&0.0006\\
					\hline
					\multirow{4}{*}{100}&2&0.0039&0.0006\\
					&3&0.011&0.0007\\
					&4&0.0361&0.0007\\
					&5&0.1215&0.0008\\
					\hline
					\multirow{4}{*}{200}&2&0.0541& 0.0009\\
					&3&0.0698& 0.0011\\
					&4& 0.1871&0.0019\\
					&5& 0.4803&0.0021\\
					\hline
					\multirow{4}{*}{500}&2&0.9866& 0.0054\\
					&3& 1.1224& 0.0070\\
					&4&1.6009&0.0087\\
					&5& 3.4348&0.0108\\
					\hline
				\end{tabular}
			}
		\end{table}
		
		To compare the accuracy of both algorithms, for each $i$th problem among 1000 problems corresponding to a fix pair $(d, cond)$, we also save the objective function value at final iteration of the two methods by $f_{\rm NESMINO}(i)$ and $f_{\rm G}(i)$, and count how many $i$ such that $|f_{\rm NESMINO}(i) - f_{\rm G}(i)|<10^{-6}$. We see that almost 1000 problems in each pair $(d, cond)$ satisfying this check.  From Table 1, we can observe that the CPU time almost  increases with $d$ and $cond$. Moreover, the smoothing algorithm depends more heavily on the shapes of ellipsoids than Gilbert's algorithm. The Table also show that both algorithm may have good potential for solving large scale problems. This example also show that despite conservative theoretical bound on the rate of convergence, in the case of ellipsoids, Gilbert's algorithm turns to be faster than smoothing algorithm .
	}
\end{Example}

\section{Conclusions}\label{sec:conr}
Minimum norm problems have been studied from both theoretical and numerical point of view in this paper.
Based on the minimum norm duality theorem, it is shown that projections onto a Minkowski sum of sets can be represented as the sum of points on constituent sets so that, at these points, all of the sets share the same normal vector. By combining Nesterov's smoothing technique and his fast gradient scheme, we have developed a numerical algorithm for solving the problems. The proposed algorithm is proved to have a better convergence rate than Gilbert's algorithm in the worst case. Numerical examples also show that the algorithm works well for the problem in high dimensions. We also note that Gilbert's algorithm is a Frank-Wolf type method; see \cite{FW56,Jag}. Although the convergence rate is known not to be
very fast, due to its very cheap computational cost per iteration, its
variants are still methods of choice in many applications.

\end{document}